\documentclass{elsart3-1}
\usepackage[utf8]{inputenc}

\usepackage{graphics}
\usepackage{graphicx}
\usepackage{epsfig}

\usepackage{amsmath,amssymb}
\usepackage[english,francais]{babel}
\RequirePackage{color}

\newtheorem{theorem}{Theorem}

\newtheorem{e-proposition}[theorem]{Proposition}

\newtheorem{e-definition}[theorem]{Definition\rm}

\renewcommand{\theequation}{\arabic{equation}}
\newcommand{\commentout}[1]{}

\renewcommand{\epsilon}{\varepsilon}
\newcommand{\eps}{\varepsilon}

\newcommand{\la}{\left\langle}
\newcommand{\ra}{\right\rangle}

\def\og{\leavevmode\raise.3ex\hbox{$\scriptscriptstyle\langle\!\langle$~}}
\def\fg{\leavevmode\raise.3ex\hbox{~$\!\scriptscriptstyle\,\rangle\!\rangle$}}

\journal{the Acad\'emie des sciences}
\begin{document}
% place in the next line the header (rubrique) chosen for your article,
% if you know it (you can also have 2, format : Header1/Header2
\centerline{}
\begin{frontmatter}


% \author{Name\thanksref{label2}}
\selectlanguage{english}
\title{Invasion fronts with variable motility:  phenotype selection, spatial sorting and wave acceleration}

% The [label1] can be suppressed if there is only one address for all authors

\selectlanguage{english}

\author[UMPA]{Emeric Bouin},
\ead{emeric.bouin@ens-lyon.fr}
\author[UMPA]{Vincent Calvez},
\ead{vincent.calvez@ens-lyon.fr}
\author[MAP5]{Nicolas Meunier},
\ead{nicolas.meunier@parisdescartes.fr}
\author[CMAP]{Sepideh Mirrahimi},
\ead{mirrahimi@cmap.polytechnique.fr}
\author[LJLL]{Beno\^it Perthame},
\ead{benoit.perthame@upmc.fr}
\author[CEFE]{Ga\"el Raoul},
\ead{g.raoul@damtp.cam.ac.uk}
\author[LPTMC]{Rapha\"el Voituriez}
\ead{voiturie@lptmc.jussieu.fr}

\address[UMPA]{Ecole Normale Sup\'erieure de Lyon, CNRS UMR 5669 UMPA, INRIA project NUMED, 
46, all\'ee d'Italie, 
F69364 Lyon}

\address[MAP5]{Universit\'e Paris Descartes, CNRS UMR 8145 MAP5, 45 rue des Saints-P\`eres, F75270 Paris} 

\address[CMAP]{Ecole Polytechnique, CNRS UMR 7641 CMAP, INRIA project MAXPLUS. Route de Saclay, F91128 Palaiseau}

\address[LJLL]{Universit\'e Pierre et Marie Curie, CNRS UMR 7598 LJLL, INRIA projet BANG, 4 Pl. Jussieu, F75005 Paris}

\address[CEFE]{Centre d'Ecologie Fonctionnelle et Evolutive, CNRS UMR 5175, 1919 route de Mende, F34293 Montpellier}  

\address[LPTMC]{Universit\'e Pierre et Marie Curie, CNRS UMR 7600 LPTMC, 4 Place Jussieu, F75252 Paris}

%\address[VC]{École Normale Sup\'erieure de Lyon,
%UMR CNRS 5669 'UMPA', and INRIA Alpes, projet NUMED,
%46, all\'ee d'Italie, 
%F-69364 LYON Cedex 07, 
%FRANCE}

% If you know the dates of reception, and acceptation you can put them now;
%  idem the name of the person presenting the Note

\medskip
\begin{center}
{\small Received *****; accepted after revision +++++\\
Presented by £££££}
\end{center}

\begin{abstract}
\selectlanguage{english}
% Text of abstract in English
%
%
%{\it To cite this article: E. Bouin, V. Calvez, C. R. Acad. Sci. Paris, Ser. I 340 (2005).}
%
%\commentout{We analyse a simple reaction-diffusion equation describing  invasion fronts in ecology,}  
Invasion fronts in ecology are well studied but very few mathematical results concern the case with variable motility (possibly due to mutations). Based on an apparently simple reaction-diffusion equation, we explain the observed phenomena of front acceleration  (when the motility is unbounded) as well as other quantitative results, such as  the selection of the most motile individuals (when the motility is bounded). The key argument for the construction and analysis of traveling fronts is the derivation of the dispersion relation linking the speed of the wave and the spatial decay. When the motility is unbounded we show that the position of the front scales as $t^{3/2}$. When the mutation rate is low we show that the canonical equation for   the dynamics of the fittest trait  should be stated as a PDE in our context. It turns out to be  a type of Burgers equation with source term. 

\vskip 0.5\baselineskip

\selectlanguage{francais}
% Text of abstract in French
\noindent{\bf R\'esum\'e} \vskip 0.5\baselineskip \noindent
{\bf  Fronts d'invasion avec motilit\'e variable: r\'epartition des ph\'enotypes et acc\'el\'eration de l'onde.} Les fronts d'invasion en \'ecologie ont \'et\'e largement \'etudi\'es. Cependant peu de r\'esultats math\'ematiques existent pour le cas d'un coefficient de motilit\'e variable (\`a cause des mutations). A partir d'un mod\`ele minimal  de r\'eaction-diffusion, nous expliquons le ph\'enom\`ene observ\'e d'acc\'el\'eration du front (lorsque la motilit\'e n'est pas born\'ee), et nous d\'emontrons l'existence d'ondes progressives ainsi que la s\'election des individus les plus motiles (lorsque la motilit\'e est born\'ee). Le point cl\'e pour la construction des fronts est  la relation de dispersion qui relie la vitesse de l'onde avec la d\'ecroissance en espace. Lorsque la motilit\'e n'est pas born\'ee nous montrons que la position du front suit une loi d'\'echelle en $t^{3/2}$. Lorsque le taux de mutation est faible, nous montrons que, dans notre contexte, l'\'equation canonique pour la dynamique du meilleur trait est une EDP. C'est une \'equation de type Burgers avec terme source. 
   
%Ce travail concerne l'\'etude de fronts de propagation pour un mod\`ele cin\'etique lin\'eaire avec un op\'erateur de relaxation de type BGK. Ce type de mod\`ele est naturellement appropri\'e dans le cadre de l'\'etude de mouvement de bact\'eries ainsi que dans le cadre de la neutronique, entres autres. On pr\'esente ici une approche de type Hamilton-Jacobi permettant de capturer la vitesse de ces fronts dans la limite hyperbolique des fronts raides. 

%{\it Pour citer cet article~: E. Bouin, V. Calvez, C. R. Acad. Sci. Paris, Ser. I 340 (2005).}

\end{abstract}
\end{frontmatter}

% now the Version franÁaise abrÃ\`agÃ\`ae, if it exists
\selectlanguage{francais}
\section*{Version fran\c{c}aise abr\'eg\'ee}

Dans cette note nous \'etudions un mod\`ele simple d\'ecrivant des fronts invasifs en \'ecologie, pour lesquels la mobilit\'e des individus est sujette \`a variations. Le mod\`ele, issu de \cite{BMV}, est le suivant,
\begin{equation}\label{eq:BMV}
\begin{cases}
\partial_t n(t,x,\theta) = \theta \partial^2_{xx} n(t,x,\theta) + r n(t,x,\theta)\left ( 1 - \rho(t,x) \right)   + \alpha \partial^2_{\theta\theta} n (t,x,\theta) \, ,  \quad x\in (-\infty, + \infty),\, \theta >0\,,
\\ 
\rho(t,x) = \int n(t,x,\theta)\, d\theta\, ,  \qquad   n(t, -\infty, \theta)=\overline N(\theta) \, , \qquad n(t, +\infty, \theta)=0 \, .
%\, , \qquad \rho(t,-\infty)=1\, .
\end{cases}
\end{equation}
Les conditions aux limites sont compl\'et\'ees ci-dessous. 
Nous \'etudions la dynamique de cette \'equation sous diff\'erents r\'egimes. Dans un premier temps nous \'etudions le probl\`eme de propagation de front \`a vitesse constante. Cela circonscrit l'espace des traits $\theta\in (0,\Theta)$, qui doit \^etre born\'e 
$\Theta<+\infty$. Nous d\'emontrons que la situation est similaire au cas de l'\'equation de Fisher-KPP. Il existe une vitesse minimale $c^*$ de propagation des ondes de r\'eaction-diffusion (Th\'eor\`eme \ref{th:sorting}). La relation de dispersion, analogue \`a celle de Fisher-KPP (un trin\^ome du second degr\'e en l'occurrence) est donn\'ee via la r\'esolution d'un probl\`eme spectral. 
%La solution s'exprime explicitement \`a l'aide des fonction d'Airy.  
On lit sur la distribution des ph\'enotypes (le vecteur propre) que les fortes motilit\'es sont favoris\'ees, conclusion oppos\'ee au cas des domaines born\'es en espace \cite{Dockery-1998}.   Ce m\^eme probl\`eme spectral intervient lorsqu'on \'etudie la propagation du front en r\'egime asymptotique hyperbolique $(t,x)\to (t/\eps,x/\eps)$ dans la limite WKB de l'optique g\'eom\'etrique, en suivant l'ansatz $n^\eps(t,x,\theta) = \exp(u^\eps(t,x)/\eps)N^\eps(t,x,\theta)$.

Dans un second temps, nous consid\'erons un espace des traits non born\'e, $\Theta = +\infty$. Dans ce cas le front acc\'el\`ere sans cesse et nous montrons heuristiquement que la loi de propagation du front est naturellement $\langle x \rangle \sim \left(\alpha^{1/4}r^{3/4}\right)t^{3/2}$. Nous exhibons une solution particuli\`ere qui confirme cette heuristique (Proposition~\ref{prop:acceleration}). 

Dans un troisi\`eme temps, nous \'etudions le r\'egime de mutations rares, et nous \'ecrivons une \'equation canonique pour l'\'evolution du {\em trait s\'electionn\'e localement} \`a l'avant du front en r\'egime asymptotique. La d\'erivation formelle de cette \'equation conduit \`a une \'equation de transport de type Burgers, avec terme source (Proposition \ref{prop:burgers}). La partie transport est d\^ue \`a la progression du front qui "transporte" les individus et donc le trait s\'electionn\'e localement. Le terme source est d\^u \`a la pression de s\'election qui tend \`a faire augmenter la valeur du trait s\'electionn\'e localement.

\selectlanguage{english}
% main text

\section*{Introduction}

Recently, several works have addressed the issue of front invasion in ecology, where the motility of individuals is subject to variability \cite{Champagnat-Meleard-2007,Arnold-Desvillettes-Prevost-2012}. It has been postulated that selection of more motile individuals can occur, even if they have no advantage regarding their reproductive rate ({\em spatial sorting}) \cite{Thomas:2004,Kokko:2006,Ronce:2007,Shine:2011}. This phenomenon has been described in the invasion of cane toads in northern Australia \cite{Phillips:2006}. It has been shown that the speed of the front increases, coincidentally with significant changes in toads morphology. Up to now, only numerical simulations have been proposed to address this issue. Here we analyse a simple model of B\'enichou et al \cite{BMV} which contains the basic features of this process: spatial mobility, logistic reproduction, and variable motility. It is given by equation \eqref{eq:BMV} where the
%\begin{equation}\label{eq:BMV}
%\partial_t n(t,x,\theta) = \theta \partial^2_{xx} n(t,x,\theta) + r n(t,x,\theta)\left ( 1 -  \rho(t,x) \right)   + \alpha \partial^2_{\theta\theta} n (t,x,\theta)\, ,  \quad x\in (-\infty, + \infty),\; \theta >0\,,
%\end{equation} where $\rho(t,x) = \int n(t,x,\theta)\, d\theta$. The
 last    term on the right hand side accounts for modifications of the dispersal rate $\theta$ of individuals due to  mutations. We consider that mutations are random and they act as a diffusion process in the phenotype space.  When needed, we impose Neumann boundary conditions in the variable $\theta$ and far-field conditions in the variable $x$.
 %: \textcolor{red}{$n(t,\pm\infty,\theta) = 0$}.
%We will consider two examples in the following: (a) {\em diffusion of the dispersal rate} $M[ n] = \partial^2_{\theta\theta} n(t,x,\theta)$ (we consider that mutations are random and they act as a diffusion process in the phenotypical space) - (b) {\em kernel distribution of the dispersal rate} $M[n] = \int G(\theta,\theta') n(t,x,\theta')\, d\theta'$, where $G$ is the phenotypical distribution of the descendant mobility, and verifies $\int G(\theta,\theta')\, d\theta = 0$. 

\section{Phenotype selection and spatial sorting in the traveling wave}

We first consider bounded dispersal rates, say $\theta\in(0,\Theta)$. Following \cite{BMV}, we seek a traveling wave solution to the equation \eqref{eq:BMV} connecting $0$ to the uniform stationary state $\overline N(\theta) \equiv \Theta^{-1}$. For $x$ large, we make the following ansatz $n(t,x,\theta) = \exp(\lambda(x - ct))Q(\theta)$, where $c>0$ is the speed of the wave, $\lambda <0$ is the spatial decay and $Q(\theta)$ denotes the phenotypic distribution of the individuals at the edge of the front. The dispersion relation is equivalent to the following spectral problem: Given a spatial decay rate $\lambda<0$, find $c(\lambda)$ and an eigenvector $Q(\theta,\lambda)$ such that 
\begin{equation}\label{eq:spectral} 
\begin{cases}
\left( \lambda c(\lambda) + \theta \lambda^2 + r\right) Q(\theta,\lambda) + \alpha \partial^2_{\theta\theta} Q(\theta,\lambda) = 0\, , \smallskip\\ 
\partial_\theta Q(0,\lambda) = \partial_\theta Q(\Theta,\lambda) = 0 \, ,\quad  \forall \theta\; Q(\theta,\lambda)\geq 0\, , \quad  \int  Q(\theta,\lambda)\, d\theta = 1\, . 
\end{cases}
\end{equation}
The  invasion front speed $c(\lambda)$  is such that $0$ is the principal eigenvalue of this spectral problem.  
Like the Fisher-KPP equation, there is a minimal speed $c^*>0$ associated with a critical spatial decay $\lambda^*<0$. 
\medskip

\begin{theorem}[Front propagation and spatial sorting]
\label{th:sorting}
For all $c\geq c^*$, there exists formally a traveling front solution $n(t,x,\theta) = N(x - ct,\theta)$. It satisfies $N(z,\theta) \sim \exp(\lambda z) Q(\theta,\lambda)$ as $z\to +\infty$.
%For a given initial data $n(0,x,\theta)$, compactly supported on $\RR\times[0,\Theta]$, the solution travels with speed $c^*$ in the following sense:
%\begin{equation}
%\begin{cases} 
%\lim_{t\to +\infty} \rho(t,x) = 1 \, , \quad \text{uniformly on $|x|<c t$}\, ,\quad \text{for $c < c^*$}\, , \\
%\lim_{t\to +\infty} \rho(t,x) = 0 \, , \quad \text{uniformly on $|x|>c t$}\, ,\quad \text{for $c > c^*$}\, . 
%\end{cases}
%\end{equation}
The phenotypic distribution at the edge of the front $Q(\theta,\lambda)$ is unbalanced towards more motile individuals. More precisely, we have 
\begin{equation}
\langle \theta \rangle_{\rm edge}(\lambda) :=  \int \theta Q(\theta,\lambda)\, d\theta > \dfrac{\Theta}2\, , 
\qquad \la \theta \ra_{\rm edge}(\lambda) \underset{\alpha \to 0}{\longrightarrow}  \Theta \,.
\label{unbal}
\end{equation}
\end{theorem}

This problem is closely related to front propagation in kinetic equations \cite{Bouin2012}. 
There is a natural extension of this result for other mutation operators as integral operators $\alpha \int G(\theta,\theta') n(t,x,\theta')\, d\theta'$. In this case, the solution of the spectral problem is deduced from the Krein-Rutman theorem.
%\cite{Bouin-preparation}.  

Interestingly, we can measure the asymmetry of the phenotypic distribution $Q(\theta,\lambda)$. The relevant quantity here is the mean diffusion coefficient at the edge of the front $\langle \theta \rangle_{\rm edge} (\lambda)$. In order to show \eqref{unbal}, we integrate the spectral problem \eqref{eq:spectral} over $(0,\Theta)$.  We get $ \lambda c + \lambda^2\langle\theta\rangle_{\rm edge}  + r  = 0$. Dividing by $Q(\theta,\lambda)$, and integrating again the spectral problem we get after integration by parts,
\begin{equation*} 
(\lambda c + r) \Theta + \frac{\Theta^2}2 \lambda^2 + \alpha \int \left|\frac{\partial_\theta Q(\theta,\lambda)}{Q(\theta,\lambda)}\right|^2\, d\theta = 0 \, .
\end{equation*}
Hence, the mean diffusion coefficient is given by
\begin{equation}\label{eq:meantheta}
\langle \theta\rangle_{\rm edge}(\lambda)= \frac{\Theta }2 + \dfrac{\alpha}{\Theta\lambda^2}\int \left|\frac{\partial_\theta Q(\theta,\lambda)}{Q(\theta,\lambda)}\right|^2\, d\theta > \frac{\Theta }2\, .
\end{equation}
Individuals at the edge of the front are more motile than at the back of the front. There, the population is homogeneous since $N(-\infty,\theta) \equiv \Theta^{-1}$ 
and the averaged diffusion coefficient is $\langle \theta\rangle_{\rm back} =  \Theta /2$. 
The estimate \eqref{eq:meantheta} measures how far the phenotypic distribution differs from the uniform distribution. 
Finally we find easily that as  $\alpha$ vanishes,  the distribution $Q(\theta,\lambda)$ concentrates as a Dirac mass on $\theta = \Theta$ hence the second statement in \eqref{unbal}: the most motile individuals are selected. The conclusion is exactly the opposite on bounded spatial domains \cite{Dockery-1998}.

We can derive additional information for the minimal speed $c^*$. Differentiating \eqref{eq:spectral} we obtain 
\[ 
\left( \lambda c'(\lambda) + c(\lambda) + 2\theta\lambda \right) Q(\theta,\lambda) + \left( \lambda c(\lambda) + \theta \lambda^2 + r\right) \frac{\partial }{\partial\lambda}Q(\theta,\lambda) + \alpha \partial^2_{\theta\theta} \frac{\partial }{\partial\lambda}Q(\theta,\lambda) = 0 \, . 
\]
Definition of $c^*$ ensures that $c'(\lambda^*) = 0$. We use the notation $\langle f \rangle:=\int f  (\theta) Q^*(\theta)\, d\theta$. Since the operator in \eqref{eq:spectral} is self-adjoint, we obtain after multiplication by $Q(\theta,\lambda)$ and integration, the relation
\begin{equation}\label{eq:cstar}
c^* \langle Q^* \rangle + 2 \lambda^*  \langle \theta Q^* \rangle = 0 \, , \quad c^* =- 2 \lambda^*  \dfrac{ \langle \theta Q^* \rangle}{\langle  Q^* \rangle}   \, ,
\end{equation}
Recalling that $ \lambda^* c^* + (\lambda^*)^2\langle\theta\rangle  + r  = 0$, we can eliminate $\lambda^*$ and we get the following expression for $c^*$,
\begin{equation*} 
(c^*)^2 =4 r \langle \theta  \rangle  \; \left( 1 - \left( 1 -  \frac{ \langle \theta  \rangle\langle Q^* \rangle }{ \langle \theta Q^* \rangle } \right)^2  \right)^{-1}  > 4 r \langle \theta  \rangle \, .
\end{equation*}
In other words, the usual formula for the  KPP wave speed  underestimates the actual minimal speed.

\section{Spatial sorting and the invasion front}

%The spectral problem \eqref{eq:spectral} also occurs when the motion of the invasion front is considered. 

Next, we focus on the invasion front. It is natural to perform the hyperbolic rescaling $(t,x)\to (t/\eps,x/\eps)$ in order to catch the motion of the front \cite{Freidlin:1986,Evans:1989}. The new equation writes after rescaling
\begin{equation*}%\label{eq:BMV-eps}
\eps \partial_t n^\eps(t,x,\theta)  =  \eps^2\theta \partial^2_{xx} n^\eps(t,x,\theta)  + r n^\eps(t,x,\theta) \left ( 1 -  \rho^\eps(t,x) \right)   + \alpha \partial^2_{\theta\theta}n^\eps(t,x,\theta)  \, .
\end{equation*}
We perform the partial WKB ansatz in the $x$ variable: $n^\eps(t,x,\theta) = \exp(u^\eps(t,x)/\eps)N^\eps(t,x,\theta)$,
with the renormalization $\int N^\eps(t,x,\theta)\, d\theta = 1$. As  $\eps\to 0$,  the first order expansion leads to solve
\begin{equation*} %\label{eq:WKB1}
\partial_t u^0(t,x)  N^0  = \theta |\partial_x u^0(t,x) |^2 N^0(t,x,\theta)  + r \left(1 -  \rho^0(t,x) \right) N^0(t,x,\theta)  + \alpha \partial^2_{\theta\theta} N^0(t,x,\theta) \, . 
\end{equation*}
The edge of the front is delimited by the area $\{u^0 (t,x)<0\}$. On this set we have $\rho^0(t,x) = 0$ by construction. Therefore we shall solve again the spectral problem \eqref{eq:spectral} for $N^0 \geq 0$.
%\begin{equation}\label{eq:eigenpb} \begin{cases}
%  \theta |p|^2 Q(p,\theta) + \alpha \partial_{\theta \theta}Q(p,\theta) = H(p) Q(p,\theta)  \, , \smallskip\\ 
%\partial_\theta Q(p,0) = \partial_\theta Q(p,\Theta) = 0 \, ,\quad \forall \theta\; Q(p,\theta)\geq 0\, , \quad  \int  Q(p,\theta)\, d\theta = 1\, . 
%\end{cases}
%\end{equation}
%where the hamiltonian $H(p)\; (p\in \RR)$ denotes the principal eigenvalue. 
Consequently the motion of the front is driven by the eikonal equation built on the effective speed $c(\lambda)$,
\begin{equation*}
\max\left( u^0, \partial_t u^0  + \partial_x u^0 \cdot c(\partial_x u^0 )  \right) = 0  \, .
\end{equation*}

The rigorous derivation of this Hamilton-Jacobi equation requires some work. We need basically refined {\em a priori} estimates on $(u^\eps)_\eps$. 
%and L. C. Evans' perturbed test function method in order to get rid of  the fast variable $\theta$. We leave it for future perspectives and refer to \cite{Bouin2012} for a similar result applied to the kinetic transport equation. 
We 
%briefly  sketch
formally show the main argument 
leading to establish  the viscosity limit $u^0$ of $u^\eps$ in the set $\{u^0 (t,x)<0\}$. We leave the complete proof for future perspectives. Let $v^0 $ be a $\mathcal C^2$ test function such that $u^0 - v^0$ has a strict maximum at $(t^0,x^0)$. 
%Let $\eta(t,x,\theta)$ be a microscopic corrector to be defined later, as in the perturbed test function of Evans. We introduce the  correction $\psi^\eps = \psi^0 + \eps \eta$. 
The function $u^\eps - v^0$ has a maximum at $(t^\eps, x^\eps)$, with $(t^\eps,x^\eps)$ close to $(t^0,x^0)$. Plugging $v^0$ into the equation satisfied by $(u^\eps,N^\eps )$, namely 
\[ \left[ \partial_t u^\eps(t,x) - \theta |\partial_x u^\eps(t,x) |^2 -  \eps\theta \partial^2_{xx} u^\eps(t,x) - r + O(\eps) \right] N^\eps(t,x,\theta)  =     \alpha \partial^2_{\theta\theta} N^\eps(t,x,\theta)   \, ,  \]
we obtain at $(t^\eps, x^\eps)$:
\begin{equation*}
\left[ \partial_t v^0(t^\eps,x^\eps) - \theta |\partial_x v^0(t^\eps,x^\eps) |^2 - r + O(\eps) \right] N^\eps(t^\eps,x^\eps,\theta)  \leq      \alpha \partial^2_{\theta\theta} N^\eps(t^\eps,x^\eps,\theta)   \, .
% \partial_t v^0(t^\eps,x^\eps) + \theta^\eps \left( |\nabla_x \psi^\eps|^2 - \eps \Delta_{xx} \psi^\eps \right) + r  + \alpha \left( \dfrac{|\partial_\theta \psi^\eps|^2}{\eps^2} - \dfrac{ \partial_{\theta\theta}^2 \psi^\eps}{\eps} \right)\, .
\end{equation*}
Therefore, $N^\eps$ is a non-negative, non-trivial subsolution of the spectral problem \eqref{eq:spectral}. From the characterization of the principal eigenvalue, we have
\[ 
\partial_t v^0(t^\eps,x^\eps) +  \partial_x v^0(t^\eps,x^\eps) \cdot c(\partial_x v^0(t^\eps,x^\eps) )  + O(\eps) \leq 0\, . 
\]
Passing to the limit   $\eps\to 0$, we obtain that $v^0$ satisfies $\partial_t v^0 +  \partial_x v^0  \cdot c(\partial_x v^0  ) \leq 0$ at $(t^0,x^0)$. 
Therefore $u^0$ is a viscosity sub-solution of the eikonal equation $\partial_t u^0  + \partial_x u^0 \cdot c(\partial_x u^0) = 0 $ in the   interior of the set $\{u^0 (t,x)<0\}$. The same argument shows that it is also a supersolution and thus a viscosity solution.

\begin{figure}
\begin{center}
\includegraphics[width = 0.46\linewidth]{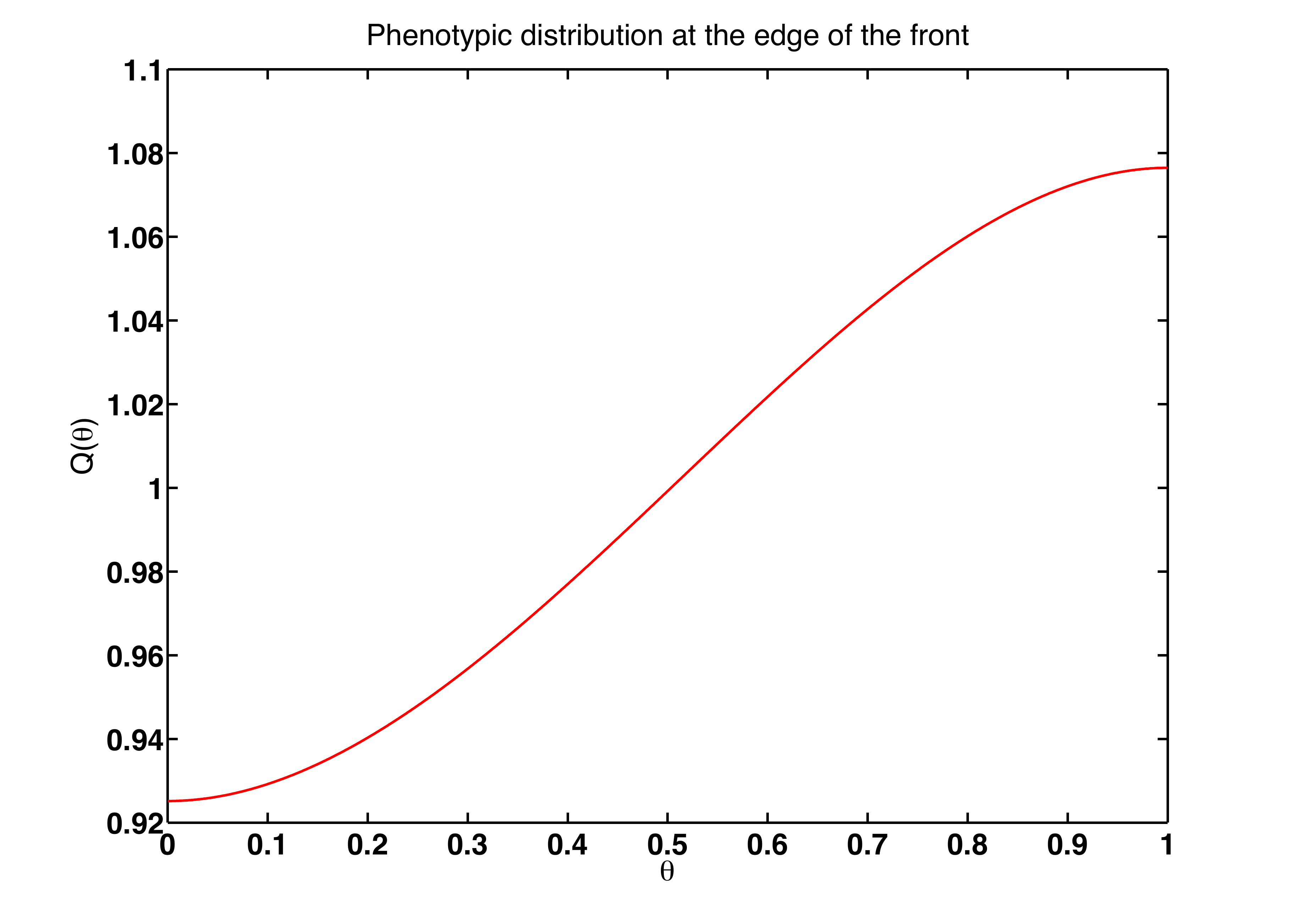}{\small (a)}
\includegraphics[width = 0.46\linewidth]{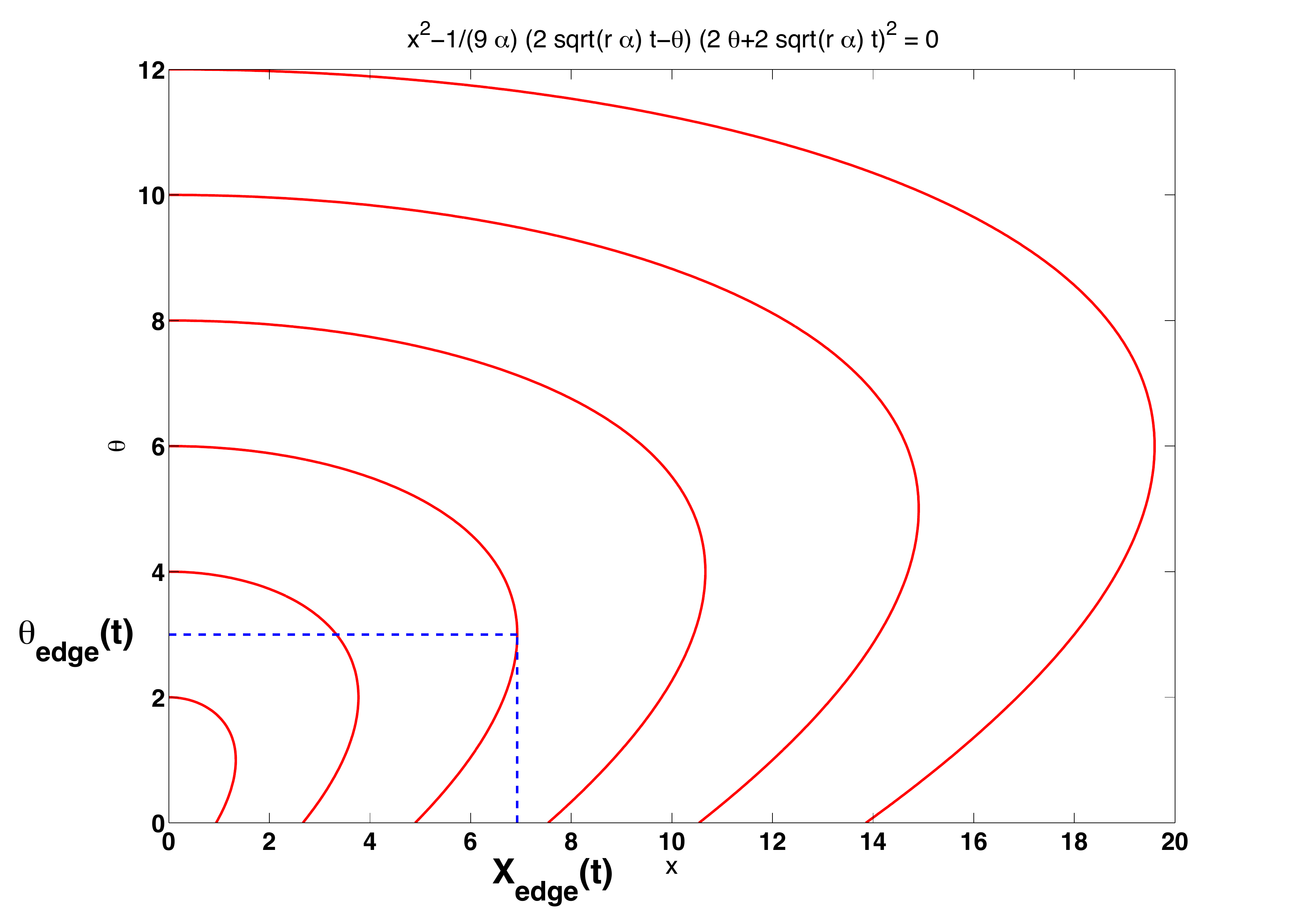}{\small (b)}
\caption{(a) {\em Phenotype selection.} At the edge of the front the phenotypic distribution is not uniform: more motile individuals are selected. We solved numerically the relation dispersion for $(r,\alpha,\Theta) = (1,1,1)$. 
%We found $c^* \approx 1.4362$ and $\lambda \approx -1.35$. 
(b) {\em Front acceleration.} The set $\{u^0  = 0\}$ is plotted in the phase space $(x,\theta)$, for successive times. The far-right point $X_{\rm edge}(t)$ determines the location of the front.}
\label{fig:nullset}
\end{center}
\end{figure}

\section{Front acceleration}

In the case where the set of dispersal rates is unbounded, say $\theta\in (0,+\infty)$, then we cannot solve the spectral problem \eqref{eq:spectral}. There is no intrinsic speed of propagation, and the front is accelerating as time goes on.
Heuristically, we expect the averaged diffusion coefficient to grow linearly with time $\langle\theta\rangle\sim \left(\sqrt{\alpha r}\right) t$  (as for the Fisher-KPP equation set in the phenotype space). Hence the invasion front should scale as $\langle x \rangle \sim \left(\alpha^{1/4}r^{3/4}\right)t^{3/2}$ since the speed is given by $ c \sim \sqrt{\langle\theta\rangle r}$. Therefore we perform the asymptotic scaling $(t,x,\theta)\to \left( t/\eps,  x/\eps^{3/2} ,  \theta/\eps \right)$, 
%and $n^\eps \to \eps n^\eps $, 
in order to catch the motion of the front asymptotically. The equation writes after rescaling
\begin{equation} \label{eq:scaling 3ter}
\eps \partial_t n^\eps  = \eps ^2\theta \partial^2_{xx} n^\eps + rn^\eps \left(1 -  \rho^\eps(t,x) \right) + \eps^2 \alpha \partial^2_{\theta\theta} n^\eps \, .
\end{equation}
We perform the WKB ansatz in both variables $(x,\theta)$: $n^\eps (t,x,\theta) = \exp\left(  u^\eps(t,x,\theta)/\eps  \right)$.
We derive formally the Hamilton-Jacobi equation for $u^0$ in the limit $\eps\to 0$: 
\begin{equation} \label{eq:double HJ}
\partial_t u^0 - \theta|\partial_x u^0|^2 - \alpha|\partial_\theta u^0|^2  =  r\left(1 -  \rho^0(t,x) \right) \, , 
\end{equation}
where $\rho^0(t,x)$ is the Lagrange multiplier associated with the constraint $ \max_\theta u^0(t,x, \theta) \leq 0$. In particular,  when the constraint is inactive, {\em i.e.} $\max_\theta u^0(t,x, \theta) <0$, we have $\rho^0(t,x)=0$. 
The behavior of the {\em phase function} $u^0(t,x,\theta)$ is as follows: the nullset $\{u^0(t,x,\theta) = 0\}$ propagates in space and phenotype under the action of space mobility, and mutations, combined with growth of individuals. We are able to compute explicitely the evolution of this set for a particular initial data, using Lagrangian Calculus. 
\medskip 

\begin{e-proposition}[Front acceleration]
\label{prop:acceleration}
For the particular initial data,
$
u^0(0,x,\theta) = \begin{cases} 0 & \text{if $(x,\theta) = (0,0)$} \\
-\infty & \text{if $(x,\theta) \neq  (0,0)$}
\end{cases}
$
the location of the nullset $\{u^0(t,x,\theta) = 0\}$ is given by the following implicit formula (see Fig.  \ref{fig:nullset})
\begin{equation*}
x^2 = \frac{1}{9 \alpha} \left( \left( 2 \sqrt{r \alpha} \right) t - \theta \right) \left( 2 \theta + \left( 2 \sqrt{r \alpha} \right) t \right)^2\, .
\end{equation*}
\end{e-proposition}

%\begin{proof}
\noindent {\em Sketch of proof.}
The Hamiltonian is given by $H((x,\theta),(p_x,p_\theta)) = \theta |p_x|^2 + \alpha |p_\theta|^2 + r$, and the corresponding Lagrangian writes 
$L((x,\theta),(v_x,v_\theta)) = v_x^2/(4\theta) + v_\theta^2/(4\alpha) - r$. 
The system of characteristics is given by $\dot X(t) = (2\theta(t) p_x(t), 2\alpha p_\theta(t))$, and $\dot P(t) = - (0, |p_x(t)|^2)$. Using Lagrangian formulation, we deduce after some calculation that, with  $Z$ the solution to equation $Z^3 + (12\theta/\alpha )Z + 24 x/\alpha =0$, $u^0(t,x,\theta)$ is given by 
%We find that $u^0$ is given by 
\[ 
u^0(t,x,\theta) = -\frac 1{4\alpha t} \left( \theta + \frac\alpha 4 Z^2 \right)^2\, .   
\]
This enables to compute the nullset $\{u^0(t,x,\theta) = 0\}$.
\qed

The far-right point of the curve is attained for $ \theta_{\rm edge} =  \left( \sqrt{r \alpha} \right) t $. This determines the location of the front.
%since $n^\eps(t,x,\theta) = \exp\left(  u^\eps(t,x,\theta)/\eps  \right)$
Hence the position of the front in space is exactly,
\begin{equation*}
X_{\rm edge}(t) =  \dfrac{4 }{ 3   }  \left(   \alpha^{1/4}   r^{3/4}  \right)  t^{3/2}\, .
\end{equation*}
%We recover the heuristic law \commentout{derived previously}   \textcolor{blue}{mentioned earlier}, $\langle x \rangle \sim \left(\alpha^{1/4}r^{3/4}\right)t^{3/2}$.

%Other scaling laws can be derived from other choices for the mutation process.

\section{Adaptive dynamics at the edge of  the front}

Interestingly enough, the same equation \eqref{eq:double HJ} can be derived in the context of adaptive dynamics, a theory which studies mutation-selection processes. It is generally assumed that the mutation process is so slow that mutants can replace the resident species before new mutants arise, if they are better adapted to their environment. This yields a {\em canonical equation} which gives the dynamical evolution of the selected trait in the population \cite{Dieckmann:1996,Champagnat:2008}. Recently, PDE-based methods have been successfully used to derive such a canonical equation for continuous mutation-selection processes \cite{Diekmann:2005,Barles:2007,Lorz-Mirrahimi}. Here we extend this theory to the case of front invasion coupled with a basic mutation process. 

The only difference with \eqref{eq:BMV} is that mutations are assumed to be rare ($\alpha\to \eps^2\alpha$):
\begin{equation*}%\label{eq:BMV-eps-2}
\partial_t n(t,x,\theta) = \theta \partial^2_{xx} n(t,x,\theta) + r n(t,x,\theta)\left ( 1 -  \rho(t,x) \right)   + \eps^2 \alpha \partial^2_{\theta\theta} n (t,x,\theta)\, .
\end{equation*}
It is natural to perform a long time rescaling $t\to t/\eps$ at the scale of evolutionary changes. Then it is useful to rescale space accordingly $x\to x/\eps$ in order to catch the motion of the front (otherwise it would travel at order $ O(1/\eps)$). With these changes of scales, we end up again with   equation \eqref{eq:scaling 3ter}, resp. \eqref{eq:double HJ} in the WKB limit.
%
%We now present an alternative viewpoint on the equation for the motion of the front \eqref{eq:double HJ}, which is closely related to adaptive dynamics. Interestingly enough, we obtain the same rescaled equation \eqref{eq:double HJ} from \eqref{eq:BMV} when the mutations are assumed to be small ($\alpha\to \eps^2\alpha$). It is natural to perform a long time asymptotics rescaling $t\to t/\eps$ at the scale of evolutionary changes. Then it is useful to rescale space accordingly $x\to x/\eps$ in order to catch the motion of the front (otherwise it would travels at order $ 1/\eps$). 
%
We restrict to the edge of the front, namely $ \sup_\theta u^0(t,x,\theta)<0$, $\rho^0(t,x)=0$. We seek a  canonical equation  for the {\em locally selected trait} $\overline{\theta}(t,x)$ such that $u^0(t,x,\overline\theta(t,x)) = \sup_\theta u^0(t,x,\theta)$.
\medskip

\begin{e-proposition}[Formal derivation of the canonical equation]
\label{prop:burgers}
The locally selected trait $\overline{\theta}(t,x)$ formally satisfies  a Burgers type equation with a source term,
\begin{equation}\label{eq:canonical burger}
  \partial_t \overline{\theta}(t,x) - 2 \left(\overline{\theta}(t,x) \partial_{x} u^0\right)  \partial_x\overline{\theta}(t,x)  = \dfrac{|\partial_x u^0|^2}{ - \partial^2_{\theta\theta} u^0 }\, .
\end{equation}
\end{e-proposition}

The speed of the transport equation is $-2 \overline{\theta}(t,x) \partial_{x} u^0$. It coincides with the local minimal speed of the traveling front, see {\em e.g. \eqref{eq:cstar}}. The positive source term accounts for the evolutionary drift which pushes the population towards higher motility (numerical simulations not shown). This equation may yield shock wave singularities, as for the classical Burgers equation, because more motile populations, when located behind less motile populations, will invade them.

{\em Proof.} 
We start from the first order condition
$\partial_\theta u^0(t,x,\overline{\theta}(t,x)) = 0$.
We differentiate this relation with respect to $t$ and $x$, respectively,
\[
\partial^2_{t\theta} u^0(t,x,\overline{\theta}(t,x)) + \left( \partial^2_{\theta \theta} u^0\right)\partial_t \overline{\theta}(t,x)  = 0\,,\qquad 
\partial^2_{x\theta} u^0(t,x,\overline{\theta}(t,x)) + \left( \partial^2_{\theta \theta} u^0\right)\partial_x \overline{\theta}(t,x)  = 0\,.
\]
On the other hand, we   differentiate equation \eqref{eq:double HJ} with respect to $\theta$,
\[ \partial^2_{\theta t} u^0 - |\partial_x u^0|^2 - 2 \theta \partial_{x} u^0 \partial^2_{\theta x} u^0 - 2 \alpha \partial_{\theta} u^0 \partial^2_{\theta \theta} u^0 = 0\, . \]
Evaluating the latter at $\theta = \overline{\theta}(t,x)$ yields
\[
\partial^2_{\theta t} u^0  - 2 \overline{\theta}(t,x) \partial_{x} u^0 \partial^2_{\theta x} u^0 =  |\partial_x u^0|^2\, .
\]
Combining these calculations, we conclude that $\overline\theta$ satisfies equation \eqref{eq:canonical burger}.
\qed
\medskip

\noindent {\em Acknowledgements.} S. M. benefits from a 2 year "Fondation Math\'ematique Jacques Hadamard" (FMJH) postdoc scholarship. She would like to thank Ecole Polytechnique for its hospitality.

\bibliography{bib-toads}
\bibliographystyle{elsart-num}
\end{document}